\documentclass[review]{elsarticle}
\usepackage{amssymb}
\usepackage{amsmath}
\usepackage{graphicx}
\usepackage{lineno,hyperref}
\usepackage{url}
\usepackage{makecell}
\usepackage{float}
\modulolinenumbers[5]
\newcommand{\tensor}[1]{\boldsymbol{\mathcal{#1}}}

\newcommand{\mat}[1]{\mathbf{#1}}
\newcommand{\vect}[1]{\mathbf{#1}}
\journal{Signal Processing: Image Communication}

\begin{document}

\begin{frontmatter}
\title{High-order Tensor Completion via Gradient-based Optimization Under Tensor Train Format}
%\title{Elsevier \LaTeX\ template\tnoteref{mytitlenote}}

\author[add1,add2]{Longhao Yuan}
\author[add2,add3]{Qibin Zhao\corref{mycorrespondingauthor}}
\author[add1]{Lihua Gui}
\author[add1,add2,add4]{Jianting Cao\corref{mycorrespondingauthor}}

\cortext[mycorrespondingauthor]{Corresponding authors}

\address[add1]{Graduate School of Engineering, Saitama Institute of Technology, Japan}
\address[add2]{Tensor Learning Unit, RIKEN Center for Advanced Intelligence Project (AIP), Japan}
\address[add3]{School of Automation, Guangdong University of Technology, China}
\address[add4]{School of Computer Science and Technology, Hangzhou Dianzi University, China}

\begin{abstract}
Tensor train (TT) decomposition has drawn people's attention due to its powerful representation ability and performance stability in high-order tensors. In this paper, we propose a novel approach to recover the missing entries of incomplete data represented by higher-order tensors. We attempt to find the low-rank TT decomposition of the incomplete data which captures the latent features of the whole data and then reconstruct the missing entries. By applying gradient descent algorithms, tensor completion problem is efficiently solved by optimization models. We propose two TT-based algorithms: Tensor Train Weighted Optimization (TT-WOPT) and Tensor Train Stochastic Gradient Descent (TT-SGD) to optimize TT decomposition factors. In addition, a method named Visual Data Tensorization (VDT) is proposed to transform visual data into higher-order tensors, resulting in the performance improvement of our algorithms. The experiments in synthetic data and visual data show high efficiency and performance of our algorithms compared to the state-of-the-art completion algorithms, especially in high-order, high missing rate, and large-scale tensor completion situations.
\end{abstract}


\end{frontmatter}

% \linenumbers

\section{Introduction}
Tensors are the high-order generalizations of vectors and matrices. Representing data by tensor can retain the high dimensional form of data and keep adjacent structure information of data. Most of the real-world data are more than two orders. For example, RGB images are order-three tensors ($height \times width \times channel$ ), videos are order-four tensors ($height \times width \times channel \times time$) and electroencephalography (EEG) signals are order-three tensors ($magnitude \times trails \times time$). When facing data with more than two orders, traditional methods usually transform data into matrices or vectors by concatenation, which leads to spatial redundancy and less efficient factorization\cite{shashua2005non}. In recent years, many theories, algorithms and applications of tensor methodologies have been studied and proposed \cite{kolda2009tensor,vasilescu2003multilinear,franz2009triplerank}. Due to the high compression ability and data representation ability of tensor decomposition, many applications related to tensor decomposition have been proposed in a variety of fields such as image and video completion \cite{acar2011scalable,zhao2015bayesian}, signal processing \cite{de2008blind,muti2005multidimensional}, brain-computer interface \cite{mocks1988topographic}, image classification \cite{shashua2001linear}, etc.

In practical situations, data missing is ubiquitous due to the error and the noise in data collecting process, resulting in the generation of data outliers and unwanted data entries. Generally, the lynchpin of tensor completion is to find the correlations between the missing entries and the observed entries. Tensor decomposition is to decompose tensor data into decomposition factors which can catch the latent features of the whole data. The basic concept of solving data completion problems by tensor decomposition is that we find the decomposition factors by the partially observed data, then we take advantages of the powerful feature representation ability of the factors to approximate the missing entries. The most studied and classical tensor decomposition models are the CANDECOMP/PARAFAC (CP)  decomposition \cite{sorber2013optimization,goulart2016tensor}, and the Tucker decomposition \cite{tucker1966some,de2000best,tsai2016tensor}. CP decomposition decomposes a tensor into a sum of rank-one tensors, and Tucker decomposition approximates a tensor by a core tensor and several factor matrices. There are many proposed tensor completion methods which employ the two tensor decomposition models. In \cite{acar2011scalable}, CP weighted optimization (CP-WOPT) is proposed. It formulates tensor completion problem as a weighted least squares (WLS) problem and uses optimization algorithms to find the optimal CP factors. Fully Bayesian CP Factorization (FBCP) in \cite{zhao2015bayesian} employs a Bayesian probabilistic model to find the optimal CP factors and CP-rank at the same time. Three algorithms based on nuclear norm minimization are proposed in \cite{liu2013tensor}, i.e., SiLRTC, FaLRTC, and HaLRTC. They extend the nuclear norm regularization for matrix completion to tensor completion by minimizing the Tucker rank of the incomplete tensor. In \cite{filipovic2015tucker}, Tucker low-$n$-rank tensor completion (TLnR) is proposed, and the experiments show better results than the traditional nuclear norm minimization methods. 

Though CP and Tucker can obtain relatively high performance in low-order tensors, due to the natural limitations of these two models, when it comes to high-order tensors, the performance of the two decomposition models will decrease rapidly. In recent years, a matrix product state (MPS) model named tensor train (TT) is proposed and becomes popular \cite{oseledets2011tensor,bengua2017efficient,yang2017tensor}. For an $N$th order tensor $\tensor{X}\in\mathbb{R}^{I_1\times \cdots \times I_N}$, CP decomposition represents data by $\tensor{O}(\sum_{n=1}^NI_nR)$ model parameters, Tucker model needs $\tensor{O}(\sum_{n=1}^NI_nR+r^N)$ model parameters, and TT model requires $\tensor{O}(\sum_{n=1}^NI_nR^2)$ parameters, where $R$ represents the rank of each decomposition model. TT decomposition scales linearly to the tensor order which is the same as CP decomposition. Though the CP model is more compact by ranks, it is difficult to find the optimal CP factors especially when the tensor order is high. Tucker model is more flexible and stable, but model parameters will grow exponentially when the tensor order increases. Tensor train is free from the `curse of dimensionality' so it is a better model to process high-order tensors. In addition to CP-based and Tucker-based tensor completion algorithms, there are several works about TT-based tensor completion. \cite{bengua2017efficient} develops the low-TT-rank algorithms for tensor completion. By tensor low-rank assumption based on TT-rank, the nuclear norm regularizations are imposed on the more balanced unfoldings of the tensor, by which the performance improvement is obtained. TT-ALS is proposed in \cite{wang2016tensor}, in which the authors employ the alternative least squares (ALS) method to find the TT decomposition factors to solve tensor completion problem. A gradient-based completion algorithm is discussed in \cite{yuan2017completion}, which is to find the TT decomposition by gradient descent method and it shows high performance in high-order tensors and high missing rates tensor completion problems. There are also tensor completion algorithms which are based on the other tensor decomposition models, i.e., tensor ring (TR) decomposition \cite{zhao2016tensor,zhao2017learning} and hierarchical Tucker (HT) decomposition. Based on TR decomposition, works in \cite{wang2017efficient,yuan2018higher,yuan2018tensor} propose algorithms named TR-ALS, TR-WOPT and TRLRF which apply ALS, gradient descent and nuclear norm minimization methods to solve various tensor completion problems. Moreover, by total variations (TV) and HT decomposition, \cite{liu2018image} proposes a completion algorithm  named STTC, which explores the global low-rank tensor structure and the local correlation structure of the data simultaneously.

In this paper, we mainly focus on developing efficient tensor completion algorithms based on TT decomposition. Though several tensor completion methods based on TT model have been proposed recently \cite{bengua2017efficient,wang2016tensor,yuan2017completion}, their applicability and effectiveness are limited. The main works of this paper are concluded as follows: 1) Based on optimization methodology and tensor train decomposition, we propose two algorithms named Tensor train Weighted Optimization (TT-WOPT) and Tensor train Stochastic Gradient Descent (TT-SGD) which apply gradient-based optimization algorithms to solve tensor completion problems. 2) We conduct simulation experiments in different tensor orders and compare our algorithms to the state-of-the-art tensor completion algorithms.  The superior performance of our algorithms is obtained in both low-order and high-order tensors. 3) We propose a tensorization method named Visual Data Tensorization (VDT) to transform visual data into higher-order tensors, by which the performance of our algorithms is improved. 4) We test the performance of our algorithms on benchmark RGB images, video data, and hyperspectral image data. The higher performance of our algorithms is shown compared to the state-of-the-art algorithms.

The rest of the paper is organized as follows. In Section 2, we state the notations applied in this paper and introduce the tensor train decomposition. In Section 3, we present the two tensor completion algorithms and analyze the computational complexities of the algorithms. In Section 4, various experiments are conducted on synthetic data and real-world data, in which the proposed algorithms are compared to the state-of-the-art algorithms. We conclude our work in Section 5.

 \section{Preliminaries and Related works}
 
 \subsection{Notations}
Notations in \cite{kolda2009tensor} are adopted in our paper. A scalar is denoted by a normal lowercase/uppercase letter, e.g., $x, X \in\mathbb{R}$, a vector is denoted by a boldface lowercase letter, e.g., $\vect{x}\in\mathbb{R}^{I}$, a matrix is denoted by a boldface capital letter, e.g., $\mat{X}\in\mathbb{R}^{I\times J}$, a tensor of order $N\geq 3$ is denoted by an Euler script letter, e.g., $\tensor{X}\in\mathbb{R}^{I_1\times I_2\times\cdots \times I_N}$. 

$\vect{x}^{(1)},\vect{x}^{(2)},\cdots,\vect{x}^{(N)}$ denotes a vector sequence, in which $\vect{x}^{(n)}$ denotes the $n$th vector in the sequence. The representations of matrix sequences and tensor sequences are denoted in the same way. An element of  tensor $\tensor{X}  \in\mathbb{R}^{I_1\times I_2\times\cdots \times I_N}$ of index $\{i_{1},i_{2},\cdots,i_{N}\}$ is denoted by $x_{i_{1}i_{2}\cdots i_{N}}$ or $\tensor{X}(i_{1},i_{2},\cdots,i_{N})$. The mode-$n$ matricization (unfolding) of tensor $\tensor{X}  \in\mathbb{R}^{I_1\times I_2\times\cdots \times I_N}$ is denoted by $\mat{X}_{(n)}\in\mathbb{R}^{I_n \times  {I_1 \cdots I_{n-1} I_{n+1} \cdots I_N}}$.

Furthermore, the inner product of two tensor $\tensor{X}$, $\tensor{Y}$ with the same size $\mathbb{R}^{I_1\times I_2\times\cdots \times I_N}$ is defined as $\langle \tensor{X},\tensor{Y} \rangle=\sum_{i_1}\sum_{i_2}\cdots\sum_{i_N}x_{i_1 i_2\cdots i_N}y_{i_1 i_2\cdots i_N}$. The Frobenius norm of $\tensor{X}$ is defined by $\left \| \tensor{X} \right \|_F=\sqrt{\langle \tensor{X},\tensor{X} \rangle}$. The Hadamard  product is denoted by `$\ast$' and it is an element-wise product of vectors, matrices or tensors of the same size. For instance, given tensors $\tensor{X}, \tensor{Y}\in\mathbb{R}^{I_1\times I_2\times\cdots \times I_N}$, $\tensor{Z}=\tensor{X}*\tensor{Y}$, then $\tensor{Z}\in\mathbb{R}^{I_1\times I_2\times\cdots \times I_N}$ and $z_{i_1 i_2 \cdots i_N}=x_{i_1 i_2 \cdots i_N} y_{i_1 i_2 \cdots i_N}$ are satisfied. The Kronecker  product of two matrices $\mat{X}\in\mathbb{R}^{I \times K}$ and $\mat{Y}\in\mathbb{R}^{J \times L}$ is $\mat{X} \otimes \mat{Y} \in\mathbb{R}^{IJ \times KL}$, see more details in \cite{kolda2009tensor}.

\subsection{Tensor Train Decomposition}
The most significant feature of TT decomposition is that the number of model parameters will not grow exponentially by the increase of the tensor order. TT decomposition is to decompose a tensor into a sequence of order-three core tensors (factor tensors): $ \tensor{G}^{(1)},\tensor{G}^{(2)},\cdots,\tensor{G}^{(N)} $. The relation between the approximated tensor $\tensor{X}  \in\mathbb{R}^{I_1\times I_2\times\cdots \times I_N}$ and core tensors can be expressed as follow:
\begin{equation}
\label{tt_decom}
\tensor{X}=\ll \tensor{G}^{(1)},\tensor{G}^{(2)},\cdots,\tensor{G}^{(N)} \gg,
\end{equation}
where for $n=1,\cdots,N$, $\tensor{G}^{(n)} \in\mathbb{R}^{R_{n-1} \times I_{n} \times R_{n}}$, $R_0=R_N=1$, and the notation $\ll \cdot \gg$ is the operation to transform the core tensors to the approximated tensor. It should be noted that, for overall expression convenience, $\tensor{G}^{(1)}\in\mathbb{R}^{ I_1 \times R_1}$ and $\tensor{G}^{(N)}\in\mathbb{R}^{R_{N-1}\times I_N}$ are considered as two order-two tensors. The sequence $R_{0}, R_{1},\cdots ,R_{N}$ is named TT-rank which limits the size of every core tensor. Furthermore, the $(i_{1},i_{2},\cdots,i_{N})$th element of tensor $\tensor{X}$ can be represented by the multiple product of the corresponding mode-$2$ slices of the core tensors as: 
\begin{equation}
\label{TT_index}
 x_{i_{1}i_{2}\cdots i_{N}}=\prod\limits_{n=1}^N{\mat{G}^{(n)}_{i_{n}}},
\end{equation} 
where $\mat{G}^{(1)}_{i_{1}}, \cdots,\mat{G}^{(N)}_{i_{N}}$ is the sequence of slices from each core tensor. For $n=1 ,2, \cdots,N$, $\mat{G}^{(n)}_{i_n}\in\mathbb{R}^{R_{n-1} \times  R_{n}}$ is the mode-$2$ slice extracted from $\tensor{G}^{(n)}$ according to each mode of the element index of $x_{i_{1}i_{2}\cdots i_N}$. $\mat{G}^{(1)}_{i_1}\in\mathbb{R}^{ R_1}$ and $\mat{G}^{(N)}_{i_{N}}\in\mathbb{R}^{R_{N-1}} $ are extracted from first core tensor and last core tensor, they are considered as two order-one matrices for overall expression convenience.

\section{Gradient-based Tensor Train Completion}

\subsection{Tensor train Weighted Optimization (TT-WOPT)}

We define $\tensor{Y}\in\mathbb{R}^{I_1\times I_2\times\cdots \times I_N}$ as the partially observed tensor with missing entries and $\tensor{X}\in\mathbb{R}^{I_1\times I_2\times\cdots \times I_N}$ is the tensor approximated by the core tensors of a TT decomposition. The missing entries of $\tensor{Y}$ are filled with zero to make $\tensor{Y}$ to be a real-valued tensor. For modeling the completion problem, the indices of the missing entries need to be specified. We define a binary tensor $\tensor{W}\in\mathbb{R}^{I_1\times I_2\times\cdots \times I_N}$ named weight tensor in which the indices of missing entries and observed entries of the incomplete tensor $\tensor{Y}$ can be recorded. Every entry of $\tensor{W}$ meets:
\begin{equation}
\label{weight}
 w_{i_{1}i_{2}\cdots i_{N}}=
 \left\{
 \begin{aligned}
 &0 \qquad \text{if} \; y_{i_{1}i_{2}\cdots i_{N}} \;\text{is a missing entry},\\
 &1 \qquad \text{if} \;  y_{i_{1}i_{2}\cdots i_{N}}\;\text{is an observed entry}.
\end{aligned}
\right.
\end{equation}
The problem of finding the decomposition factors of an incomplete tensor can be formulated by a weight least squares (WLS) model. Define $\tensor{Y}_w=\tensor{W}\ast \tensor{Y}$, and $\tensor{X}_w=\tensor{W}\ast \tensor{X}$, then the WLS model for calculating tensor decomposition factors is formulated by:
\begin{equation}
\label{of1}
f(\tensor{G}^{(1)},\tensor{G}^{(2)},\cdots,\tensor{G}^{(N)})=\frac{1}{2} \left \|\tensor{Y}_w-\tensor{X}_w \right \|^{2}_F.
\end{equation}
This is an optimization objective function w.r.t. all the TT core tensors and we aim to solve the model by gradient descent methods. The relation between the approximated tensor $\tensor{X}$ and the TT core tensors can be deduced as the following equation \cite{cichocki2016tensor}:
\begin{equation}
\label{XG_relation}
\mat{X}_{(n)}=\mat{G}_{(2)}^{(n)}(\mat{G}_{(1)}^{>n} \otimes \mat{G}_{(n)}^{<n}),
\end{equation}
where for $n=1,...,N$,
\begin{equation}
\mat{G}^{>n}=\ll \tensor{G}^{(n+1)},\tensor{G}^{(n+2)},\cdots,\tensor{G}^{(N)} \gg \in\mathbb{R}^{R_n\times I_{n+1}\times \cdots \times I_N},
\end{equation}
\vspace{-0.8cm}
\begin{equation}
\mat{G}^{<n}=\ll \tensor{G}^{(1)},\tensor{G}^{(2)},\cdots,\tensor{G}^{(n-1)} \gg \in\mathbb{R}^{I_1 \times \cdots  \times I_{n-1} \times R_{n-1}}.
\end{equation}
$\mat{G}^{>n}$ and $\mat{G}^{<n}$ are the tensors generated by merging the selected TT core tensors, and we define $\mat{G}^{>N}=\mat{G}^{<1}=1$.

By equation (\ref{XG_relation}), for $n=1,...,N$, the partial derivatives of the objective function (\ref{of1}) w.r.t. the mode-2 matricization of the $n$th core tensor $\tensor{G}^{(n)}$ can be inferred as:
\begin{equation}
\label{derTT-WOPT}
\frac{\partial{f}}{\partial{\mat{G}_{(2)}^{(n)}}}=(\mat{X}_{w(n)}-\mat{Y}_{w(n)})(\mat{G}_{(1)}^{>n} \otimes \mat{G}_{(n)}^{<n})^ \mathrm{ T }.
\end{equation}

After the objective function and the gradients are obtained, we can apply various optimization algorithms to optimize the core tensors. The implementation procedure of TT-WOPT to find the TT decomposition from incomplete tensor $\tensor{Y}$ is listed in Algorithm 1.

\begin{table}[h]
\footnotesize
\begin{center}
\begin{tabular}{l}
\hline
\textbf{Algorithm 1} Tensor train Weighted Optimization (TT-WOPT)\\
\hline
1: \;\;  \textbf{Input}: incomplete tensor $\tensor{Y}$, weight tensor $\tensor{W}$ and TT-rank $\vect{r}$.\\
2: \;\;  Randomly initialize the core tensors $\tensor{G}^{(1)},\tensor{G}^{(2)},\cdots,\tensor{G}^{(N)} $.\\
3: \;\; \textbf{While} the optimization stopping condition is not satisfied \\
4: \;\; Compute $\tensor{X}_w=\tensor{W}\ast \ll \tensor{G}^{(1)},\tensor{G}^{(2)},\cdots,\tensor{G}^{(N)}\gg$.\\
5:  \;\; \textbf{For} n=1:N \\
6:  \;\; \; Compute gradients according to equation (\ref{derTT-WOPT}).\\
7: \;\; \textbf{End}\\
8: \;\; Update $\tensor{G}^{(1)},\tensor{G}^{(2)},\cdots,\tensor{G}^{(N)} $ by gradient descend method. \\
9: \;\; \textbf{End while}\\
10: \ \textbf{Output}: $\tensor{G}^{(1)},\tensor{G}^{(2)},\cdots,\tensor{G}^{(N)} $.\\
\hline
\end{tabular}
\end{center}
\end{table}

\subsection{Tensor Train Stochastic Gradient Descent (TT-SGD)}

As seen from equation (\ref{of1}), TT-WOPT computes the gradients by the whole scale of the tensor for every iteration. The computation can be redundant because the missing entries still occupy the computational space. If the scale of data is huge and the number of missing entries is high, then we only need to apply a small amount of the observed entries. In this situation, TT-WOPT can waste much computational storage and the computation will become time-consuming. In order to solve the problems of TT-WOPT as mentioned above, we propose the TT-SGD algorithm which only randomly samples one observed entry to compute the gradients for every iteration.

Stochastic Gradient Descent (SGD) has been applied in matrix and tensor decompositions \cite{gemulla2011large,maehara2016expected,wang2016online}. For every optimization iteration, we only use one entry which is randomly sampled from the observed entries, and one entry can only influence the gradient of part of the core tensors. For one observed entry of index $\{i_1, i_2, \cdots i_N \}$, if a value approximated by TT core tensors is $x_{i_1i_2 \cdots i_N}$ and the observed value (real value) is $y_{i_1i_2\cdots i_N}$, by considering equation (\ref{TT_index}), the objective function can be formulated by: 
\begin{equation}
\label{so}
f(\mat{G}^{(1)}_{i_1}, \mat{G}^{(2)}_{i_2}, \cdots, \mat{G}^{(N)}_{i_N})=\frac{1}{2}\left\|y_{i_1i_2 \cdots i_N}-\prod\limits_{k=1}^N{\mat{G}^{(k)}_{i_{k}}}\right\|^2_F.
\end{equation}
For $n=1,2,\cdots,N$, the partial derivatives of every corresponding slice $\mat{G}^{(n)}_{i_n}$ w.r.t. index $\{i_1, i_2, \cdots i_N \}$ is calculated as:
\begin{equation}
\label{derTT-SGD}
\frac{\partial{f}}{\partial{\mat{G}^{(n)}_{i_n} }}=(x_{i_1i_2 \cdots i_N}-y_{i_1i_2 \cdots i_N})(\prod\limits_{k=n+1}^N{\mat{G}^{(k)}_{i_{k}}}\prod\limits_{k=1}^{n-1}{\mat{G}^{(k)}_{i_{k}}})^T.
\end{equation}
From the equation we can see, the computational complexity of TT-SGD is not related to the scale of the observed tensor or the number of observed entries, so it can process large-scale data by much smaller computational complexity than TT-WOPT. This algorithm is also suitable for online/real-time learning. The optimization process of TT-SGD is listed in Algorithm 2:
\begin{table}[h]
\footnotesize
\begin{center}
\begin{tabular}{l}
\hline
\textbf{Algorithm 2} Tensor Train Stochastic Gradient Descent (TT-SGD)\\
\hline
1:\;\; \textbf{Input}: incomplete tensor $\tensor{Y}$ and $TT-rank$ $\vect{r}$.\\
2:\;\; Randomly initialize core tensors $\tensor{G}^{(1)},\tensor{G}^{(2)},\cdots,\tensor{G}^{(N)} $.\\
3:\;\; \textbf{While} the optimization stopping condition is not satisfied \\
4:\;\; Randomly sample $y_{i_1i_2\ldots i_N}$ from $\tensor{Y}$.\\
5:\;\; \textbf{For} n=1:N \\
6: \;\;\; Compute the gradients of the core tensors by equation (\ref{derTT-SGD}).
 \\
7:\;\; \textbf{End} \\
8:\;\; Update $\mat{G}^{(1)}_{i_1}, \mat{G}^{(2)}_{i_2}, \cdots, \mat{G}^{(N)}_{i_N}$ by gradient descent method. \\
9:\;\; \textbf{End while} \\
10:  \textbf{Output}: $\tensor{G}^{(1)},\tensor{G}^{(2)},\cdots,\tensor{G}^{(N)} $.\\
\hline
\end{tabular}
\end{center}
\end{table}

\subsection{Computational Complexity}
For tensor $\tensor{X}  \in\mathbb{R}^{I_1\times I_2\times\cdots \times I_N}$, we assume all $I_1,  I_2, \cdots, I_N$ is equal to $I$, and $R_1=R_2= \cdots =R_{N-1}=R$. According to equation (\ref{derTT-WOPT}) and (\ref{derTT-SGD}), the time complexity of TT-WOPT and TT-SGD are $\tensor{O}(NI^N+NI^{N-1}R^2)$ and $\tensor{O}(N^2R^3)$ respectively, and the space complexity of the two algorithms is $\tensor{O}(I^N+I^{N-1}R^2)$ and $\tensor{O}(R^2)$ respectively. Though TT-WOPT has larger computational complexity, it has a steady and fast convergence when processing normal-size data. TT-SGD is free from data dimensionality and the complexity of every iteration is extremely low, so it is more suitable to process large-scale data. It should be noted that for every iteration of TT-SGD, we can also apply the batch-based SGD method which calculates the summation of the gradients of bath-sized entries for every iteration. Though this can improve the stability of TT-SGD and the algorithm might need fewer iterations to be converged, the computational complexity will be increased and more computational time is needed for every iteration. In this paper, we only apply batch-one SGD algorithm, and the synthetic experiment in the next section show that our method can also achieve fast and stable convergence. The code of the proposed algorithms is available at $https://github.com/yuanlonghao/T3C\_tensor\_completion$.

\section{Experiment results}
In this section, simulation experiments are conducted to show the performance of our algorithms and the compared algorithms under various tensor orders. For real-world data experiments, we test our algorithms by color images, video data and hyperspectral image data. TT-WOPT and TT-SGD are compared with several state-of-the-art algorithms: TT-ALS \cite{wang2016tensor}, SiLRTC-TT \cite{bengua2017efficient}, TRALS \cite{wang2017efficient}, STTC \cite{liu2018image}, CP-WOPT \cite{acar2011scalable}, FBCP \cite{zhao2015bayesian}, HaLRTC and FaLRTC \cite{liu2013tensor}, and TLnR \cite{filipovic2015tucker}. For all the compared algorithms, the input incomplete tensor is $\tensor{W}*\tensor{Y}$, where $\tensor{Y}$ is the fully observed true tensor, $\tensor{W}$ is the binary tensor recording the position of observed entries. The final completed tensor $\tensor{Z}$ is calculated by $\tensor{Z}=(1-\tensor{W})*\tensor{X}+\tensor{W}*\tensor{Y}$, where $\tensor{X}$ is the output tensor obtained by each algorithm. We apply relative squared error (RSE) which is defined as $RSE=\left \| \tensor{Y}-\tensor{Z} \right \|_F/\left \| \tensor{Y} \right \|_F$ to evaluate the completion performance for each algorithm. For experiments of random missing cases, we randomly remove data points according to different missing rates $m_r$ which is defined as $m_r=1-M/\prod_{n=1}^{N}{I_n}$, where $M$ is the number of the observed entries. Moreover, to evaluate the completion quality of visual data, we introduce PSNR (Peak Signal-to-noise Ratio). PSNR is obtained by $PSNR=10\log_{10}(255^2/MSE)$, where MSE is deduced by $MSE=\Vert \tensor{Z}-\tensor{Y} \Vert_F^2/num(\tensor{Z})$, and $num(\cdot)$ denotes the number of the element of the tensor. 

For optimization method of TT-WOPT, in order to have a clear comparison with CP-WOPT which is also based on gradient descent methods, we adopt the same optimization method as paper \cite{acar2011scalable}. The paper applies nonlinear conjugate gradient (NCG) with Hestenes-Stiefel updates \cite{wright1999numerical} and the Mor\'e-Thuente line search method \cite{more1994line}. The optimization method is implemented by an optimization toolbox named Pablano Toolbox \cite{dunlavy2010poblano}. For TT-SGD, we employ an algorithm named Adaptive Moment Estimation (Adam) as our gradient descent method, it has prominent performance on stochastic-gradient-based optimization   \cite{kingma2014adam,ruder2016overview}. The update rule of Adam is as follow:
 \begin{equation}
\theta_{t+1}=\theta_t-\frac{\eta}{\sqrt{v_t}+\epsilon}m_t,
 \end{equation}
where $t$ is the iteration time of optimization value $\theta$, $\eta$ and $\epsilon$ are hyper parameters, $m_t$ and $v_t$ are the first moment estimate and second moment estimate of gradient $g_t$ respectively. $m_t=\beta_1m_{t-1}+(1-\beta_1)g_t$, $v_t=\beta_2v_{t-1}+(1-\beta_2)g_t^2$, where $\beta_1$ and $\beta_2$ are hyper parameters. For choosing the hyper parameters in Adam method, we adopt the reference values from paper \cite{kingma2014adam}. The values of $\beta_1$, $\beta_2$ and $\epsilon$ are set as 0.9, 0.999 and $10^{-8}$ respectively. The selection of learning rate is essential to the convergence speed and the performance of the gradient-based algorithms,  in our experiments, we empirically choose the learning rate $\eta$ from $\{0.0001, 0.0005, 0.001\}$ to obtain the best convergence speed and the best performance. In addition, all the data in our experiments are regularized to 0 to 1 to make the algorithms more effective.

We mainly adopt two optimization stopping conditions for all the compared completion algorithms. One is the error of two adjacent iterations of the objective function value: $|f_t-f_{t-1}| \leq tol$, where $f_t$ is the objective function value of the $t$th iteration and we set $tol=1e-4$ in our experiment. The other stopping condition is the maximum number of iteration which is set according to the scale of data and different algorithms, e.g., the maximum iteration for most algorithms are set as $500$ and for TT-SGD it usually set from $10^5$ to $10^7$. If one of the two conditions is satisfied, the optimization will be stopped. All the computations are conducted on a Mac PC with Intel Core i7 and 16GB DDR3 memory, and the computational time of the algorithms are recorded in some experiments based on this configuration.

\subsection{Synthetic Data}
We apply synthetic data generated from a highly oscillating function: $f(x)=sin\frac{x}{4}cos(x^2)$ \cite{khoromskij2015tensor} in our simulation experiments. The synthetic data is expected to be well approximated by tensor decomposition models. We sample $I^N$ entries from the values generated from the function, then the sampled values are reshaped to the desired tensor size. We employ four different tensor structures: $26 \times 26 \times 26$ (3D), $7 \times7 \times7 \times7 \times7$ (5D), $4 \times4 \times4 \times4 \times4 \times4 \times4$ (7D), and $3 \times3 \times3 \times3 \times3 \times3 \times3 \times3$ (9D), then we test TT-SGD ,TT-WOPT, TT-ALS, SiLRTC-TT, TR-ALS, CP-WOPT and HaLRTC on the synthetic data. For parameter settings, the hyper-parameters of each algorithm are tuned to obtain the best performance. For simplicity, we set values of each TT-rank and TR-rank identically, i.e., $R_1=\cdots =R_{N-1}$ for TT and $R_1=\cdots =R_N$ for TR. Moreover, the TT-rank, TR-rank and CP-rank are set as $12,10$ and $30$ under all the different tensor orders for the corresponding algorithms to make a clear comparison of the completion performance. In addition, the maximum iteration of TT-SGD is set as $10^5$, and iteration for other algorithms are all set as $500$.

\begin{figure}[h]
\begin{center}
\includegraphics[width=1\linewidth]{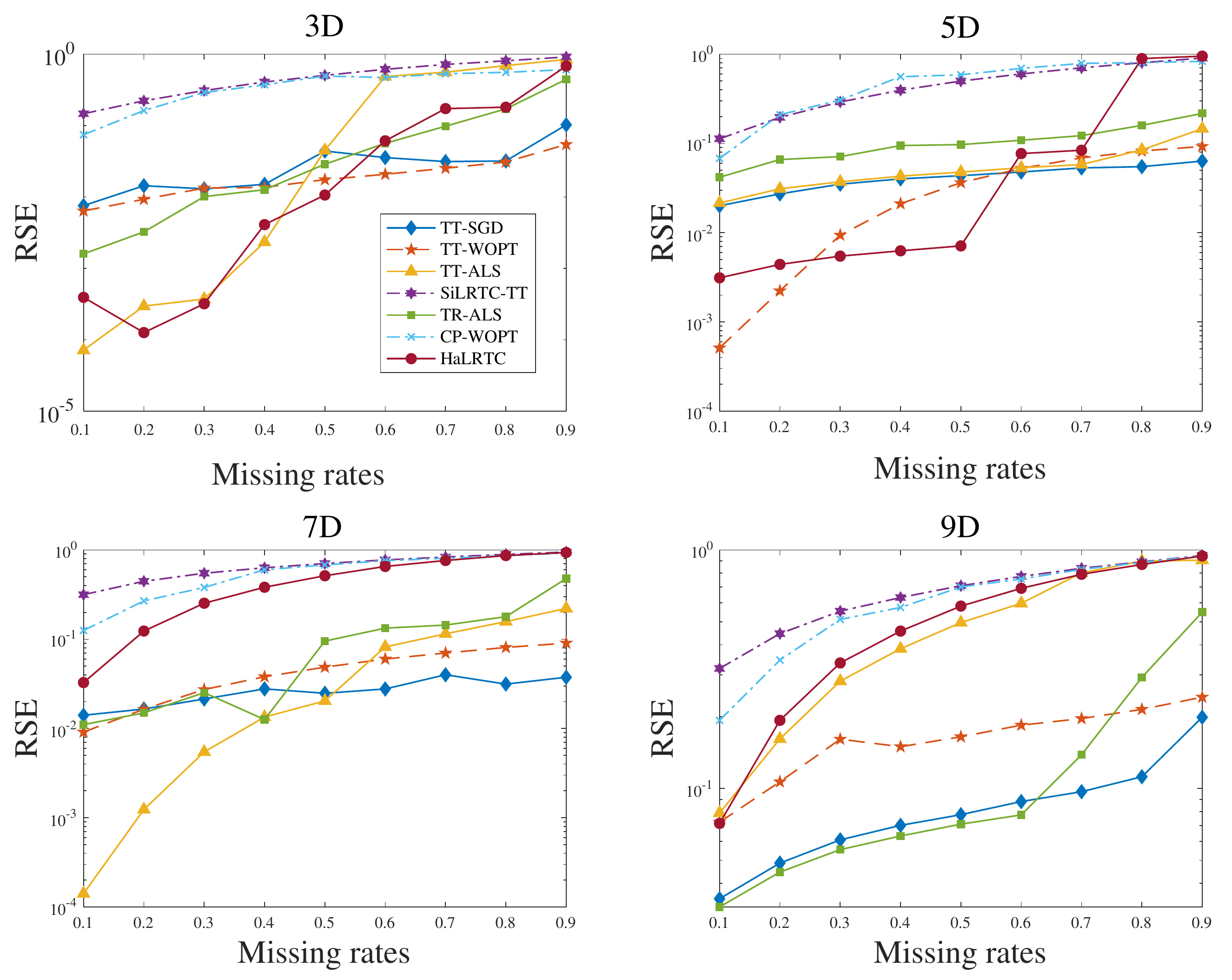}
\caption{RSE comparison of seven algorithms under four different tensor orders. The missing rate is tested from 0.1 to 0.9.}
\label{Sim}
\end{center}
\end{figure}

The graphs of Figure \ref{Sim} show the experiment results of RSE values, which change by different $m_r$ (from 0.1 to 0.9) under the four different tensor orders. From the figure, we can see that TT-WOPT and TT-SGD show high performance in all the cases. HaLRTC only shows high performance in 3D tensor case, and CP-WOPT and SiLRTC show stable but low performance in every case. Though TT-ALS and TR-ALS show higher performance than our algorithms in some low missing rate cases, the drastic performance decrease can be obtained from them when the missing rate increases, and our algorithms always show high and stable performance. 

\begin{figure}[h]
\begin{center}
\includegraphics[width=1\linewidth]{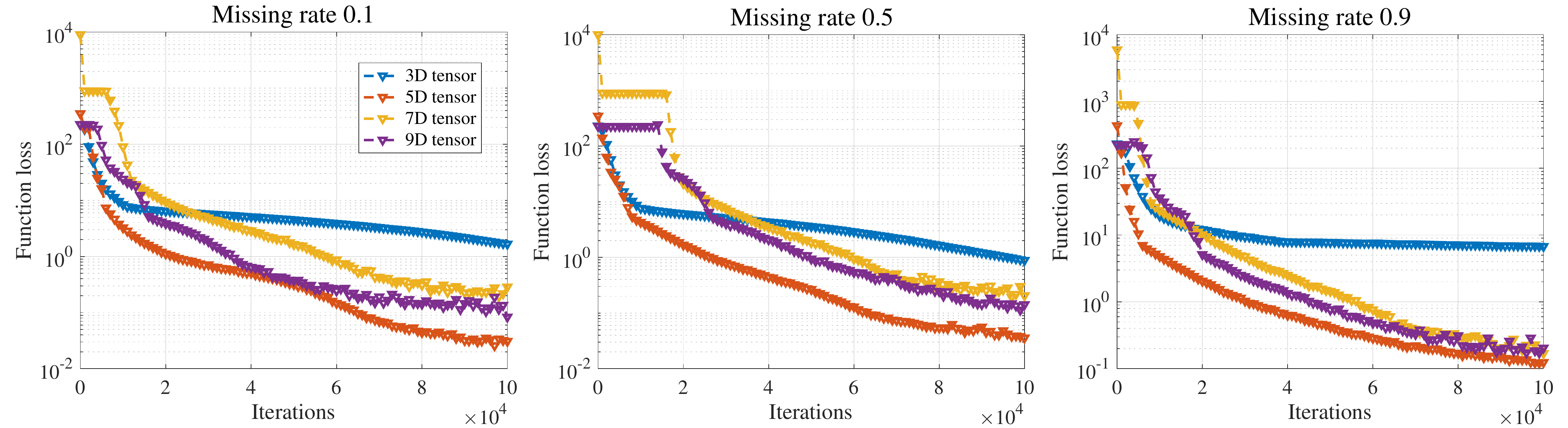}
\caption{Convergence performance of TT-SGD under four different synthetic tensors. From left to right, the missing rate of the data in each figure is 0.1, 0.5 and 0.9 respectively.}
\label{sgd}
\end{center}
\end{figure}

For the next synthetic data experiment, we aim to look into the convergence performance of the proposed TT-SGD. The four tensors which applied in the previous experiment is employed as the input data. We record the value of loss function (i.e., $\frac{1}{2}\Vert \tensor{Z}-\tensor{Y}\Vert_F^2 $) for every $10^3$ iterations and Figure \ref{sgd} shows the convergence status of TT-SGD when the missing rate is 0.1, 0.5 and 0.9 respectively. Though our TT-SGD needs large numbers of iteration to be converged, the computational complexity of each iteration is rather low (i.e., $N^2R^3$), and only one entry is sampled to calculate the gradient for every iteration. For TT-SGD, the running time of reaching $10^5$ iterations for the 3D, 5D, 7D, 9D data under the parameter setting in the experiment is 10.09 seconds, 25.09 seconds, 45.86 seconds and 75.41 seconds respectively, while for TT-WOPT, it takes about two times longer than TT-SGD (i.e., 18.80 seconds, 41.84 seconds, 100.02 seconds and 122.77 seconds) to converge to the same RSE values. The performance and computation time manifest the effectiveness of the TT-SGD algorithm.

\subsection{Visual Data Tensorization (VDT) method}

\begin{figure}[h]
\begin{center}
\includegraphics[width=0.9\linewidth]{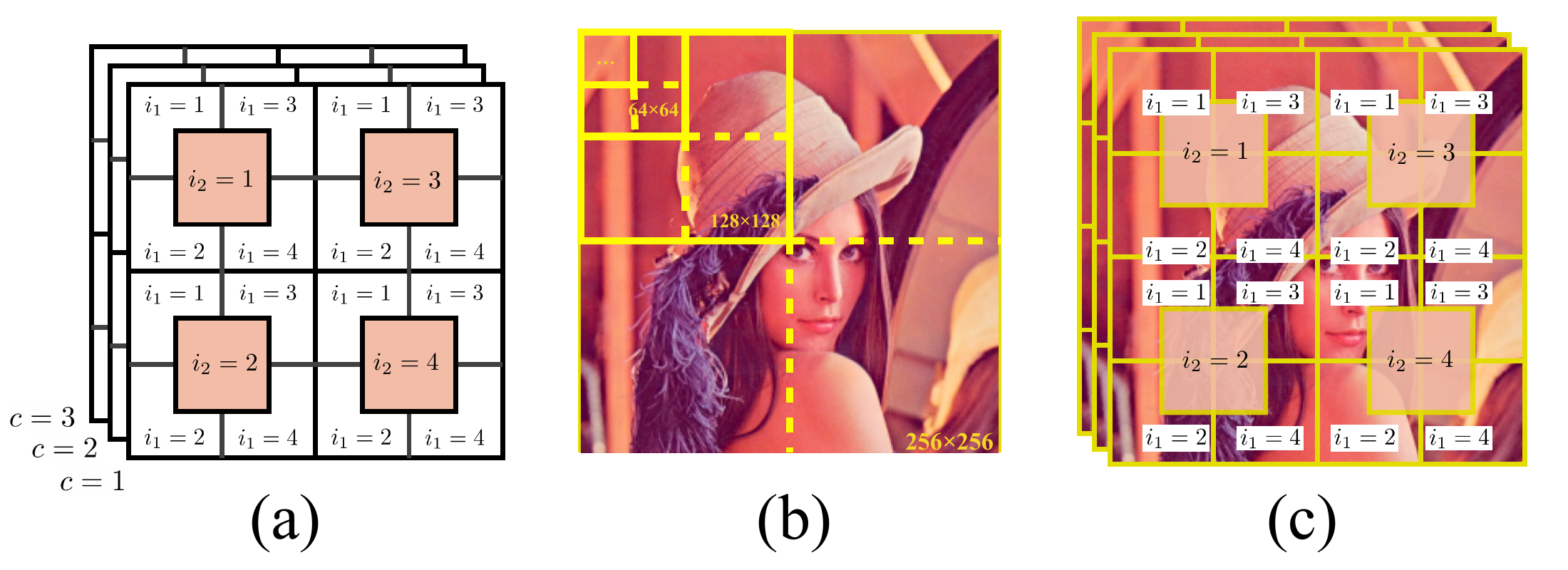}
\caption{Illustration of the proposed VDT method. Figure (a) is the example of applying VDT method on an $I\times I \times C$ tensor. Figure (b) and Figure (c) shows the example of the VDT operation on a $256 \times 256 \times 3$ image.}
\label{VDT}
\end{center}
\end{figure}

From the simulation results we can see, our proposed algorithms achieve high and stable performance in high-order tensors. In this section, we provide a Visual Data Tensorization (VDT) method to transform low-order tensor into higher-order tensor and improve the performance of our algorithms. The VDT method is derived from an image compression and entanglement methodology \cite{latorre2005image} which is to transform a gray-scale image of size $2^l \times 2^l$ into a real ket of a Hilbert space. The method cast the image to a higher-order tensor structure with an appropriate block structured addressing. Similar method named KA augmentation is proposed in \cite{bengua2017efficient} which extends the method in \cite{latorre2005image} to order-three visual data of size $2^l \times 2^l \times 3$. Our VDT method is a generalization of the KA augmentation,  and the visual data of various data sizes can be applied to our tensorization method. For visual data like RGB image, video, hyperspectral image, the first two orders of the tensor (e.g., $\mat{Y} \in\mathbb{R}^{U\times V }$) are named as the image modes. The 2D representation of the image modes cannot fully exploit the correlation and local structure of the data, so we propose the VDT method to strengthen the local structure correlation of visual data. The VDT method operates as follows: if the first two orders of a visual data tensor is $U\times V$ and can be reshaped to $u_1 \times u_2 \times \cdots \times u_l \times v_1 \times v_2 \times \cdots \times v_l$, then VDT method permutes and reshapes the data to size $u_1v_1 \times u_2v_2 \times \cdots \times u_l v_l $ and obtain the higher-order representation of the visual data. This higher-order tensor is a new structure of the original data: the first order of this higher-order tensor corresponds to a $u_1 \times v_1$ pixel block of the image, and the following orders of $u_2v_2, \cdots, u_{l}v_{l}$ describe the expanding larger-scale partition of the image. Based on VDT method, TT-based algorithms can efficiently exploit the structure information of visual data and achieve a better low-rank representation. After the tensorized data is calculated by the completion algorithms, a reverse operation of VDT is conducted to get the original image structure. The diagrams to explain the procedure of VDT are shown in Figure \ref{VDT}.

\begin{figure}[ht]
\begin{center}
\includegraphics[width=0.8\linewidth]{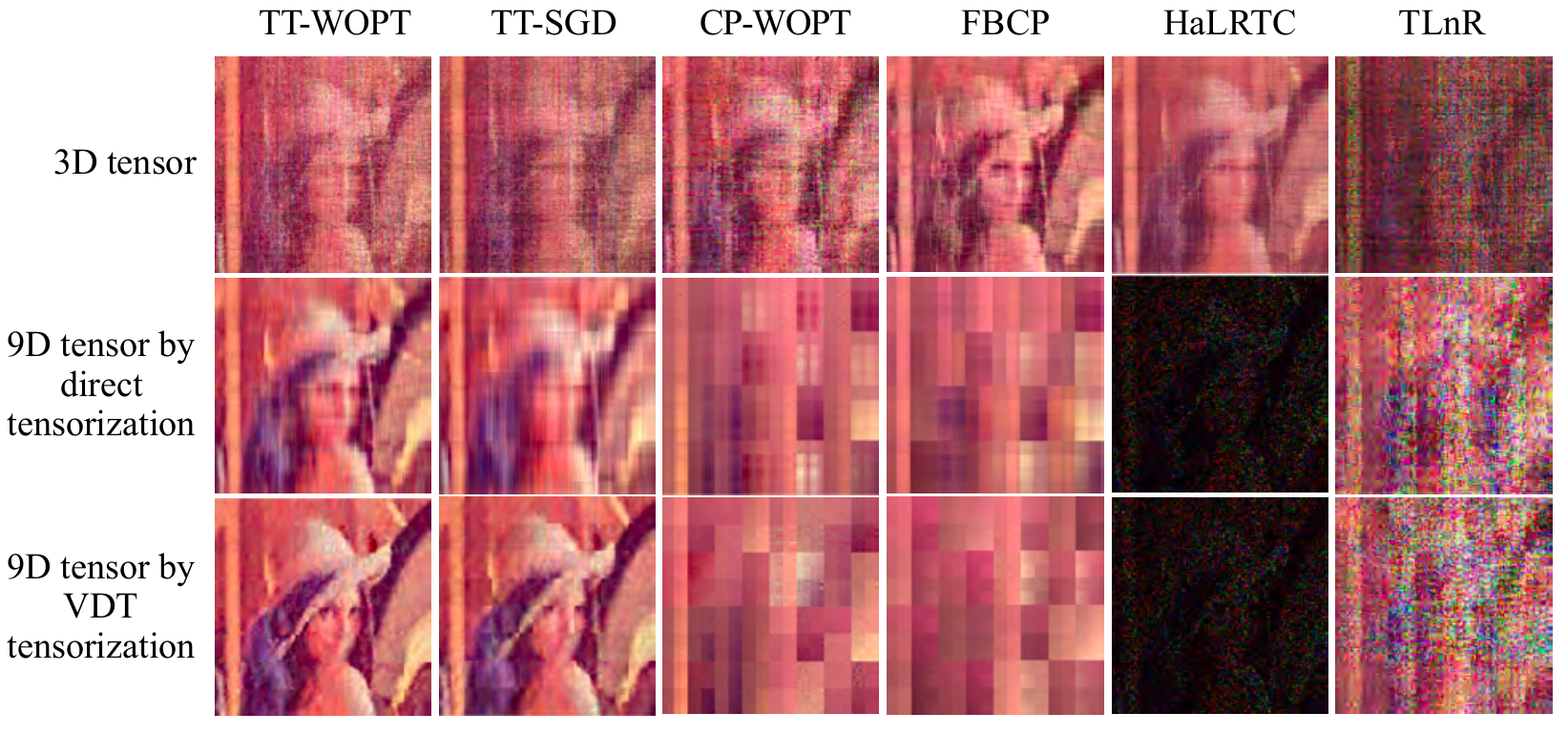}
\caption{Visual results for completion of the 0.9 random missing `Lena' image under six algorithms. The first row applies original order-three tensor data, the second row applies order-nine tensor data without VDT method, and the third row applies order-nine tensor data generated by VDT method.}
\label{90}
\end{center}
\end{figure}

To verify the effectiveness of our VDT method, we choose a benchmark image `Lena' with $0.9$ missing rate. We compare the performance of the six algorithms (TT-WOPT, TT-SGD, CP-WOPT, FBCP, HaLRTC and TLnR) under three different data structures: order-three tensor, order-nine tensor without VDT, order-nine tensor generated by VDT method. The order-three tensor applies original image data structure of size $256\times 256\times 3$. The nine-order tensor without VDT is generated by directly reshaping data to the size $4 \times4 \times4 \times4 \times4 \times4 \times4 \times4 \times3$. For nine-order tensor with VDT method, firstly the original data is reshaped to a order-seventeen tensor of size $2 \times2 \times 2 \times2 \times2 \times2 \times2 \times2 \times2 \times2 \times2 \times2 \times2 \times2 \times2 \times2 \times3 $ and then it is permuted according to the order of $\{1 \;9 \;2 \;10 \;3 \;11 \;4 \;12 \;5 \;13 \;6 \;14 \;7 \;15 \;8 \;16 \;17\}$. Finally we reshape the tensor to a nine-order tensor of size $4 \times4 \times4 \times4 \times4 \times4 \times4 \times4 \times3$. This nine-order tensor with VDT is considered to be a better structure of the image data. The first order of the nine way tensor contains the data of a $2 \times 2$ pixel block of the image and the following orders of the tensor describe the expanding pixel blocks of the image. Most of the parameter settings follow the previous synthetic data experiments, and we tune the TT-rank, CP-rank and Tucker-rank of the corresponding algorithms to obtain the best performance. Figure \ref{90} and Table \ref{90t} show the visual results and numerical results of the six algorithms under the three different data structure. We can see that in the three-order tensor case, the results among the algorithms are similar. However, for nine-order cases, other algorithms fail the completion task while TT-WOPT and TT-SGD perform well. Furthermore, when the image is transformed to nine-order tensor by VDT method, we see the distinct improvement of our two algorithms. 

\begin{table}[h]
\centering
\caption{Numerical results of completion performance (RSE and PSNR) of six algorithms under three tensor structures of image `Lena'. }
\label{90t}
\resizebox{\textwidth}{20mm}{
\begin{tabular}{c|c|c|c|c|c|c|c}
\hline
\hline
 \multicolumn{2}{c|}{}& TT-WOPT & TT-SGD & CP-WOPT & FBCP & HaLRTC&TLnR\\
 \hline
three-order&\makecell[cc]{RSE \\ PSNR}   &     \makecell[cc]{0.2822\\ 16.12}    &    \makecell[cc]{0.2604\\ 16.84}     &
\makecell[cc]{0.3392\\ 14.53}     &           \makecell[cc]{\textbf{0.1942}\\ \textbf{19.36}}   &\makecell[cc]{0.1981\\  19.18}&\makecell[cc]{0.6552\\ 8.802}\\
\hline
nine-order&\makecell[cc]{RSE \\ PSNR}&\makecell[cc]{\textbf{0.1558}\\ \textbf{21.31}}&\makecell[cc]{0.1793\\ 20.06}&\makecell[cc]{0.2562\\ 16.95}&\makecell[cc]{0.2682\\ 16.57} &   \makecell[cc]{0.9310\\ 5.746}&\makecell[cc]{1.207\\ 3.486}\\
\hline
nine-order VDT&\makecell[cc]{RSE \\ PSNR}&\makecell[cc]{\textbf{0.1262}\\ \textbf{23.21}}&\makecell[cc]{0.1493 \\ 21.77}&\makecell[cc]{0.2573\\16.97}&\makecell[cc]{0.2687\\ 16.57}  & \makecell[cc]{0.9301\\ 5.751}&\makecell[cc]{0.7114\\ 10.84}\\
\hline\hline
\end{tabular}
}
\end{table}

\subsection{Benchmark Image Completion}

From the previous experiments we can see, TT-based and TR-based algorithms can be applied to higher-order tensors, and significant improvement of TT-based algorithms can be seen when the VDT method is applied to the image tensorization. However, for algorithms which are based on CP decomposition and Tucker decomposition, higher-order tensorization will decrease the performance. In later experiments, we only apply the VDT method to TT-WOPT, TT-SGD, TT-ALS, SILRTC-TT and TR-ALS. For CP-WOPT, FBCP, TLnR, STTC and HaLRTC, we keep the original data structure to get better results.

\begin{figure}[H]
\begin{center}
\includegraphics[width=1\linewidth]{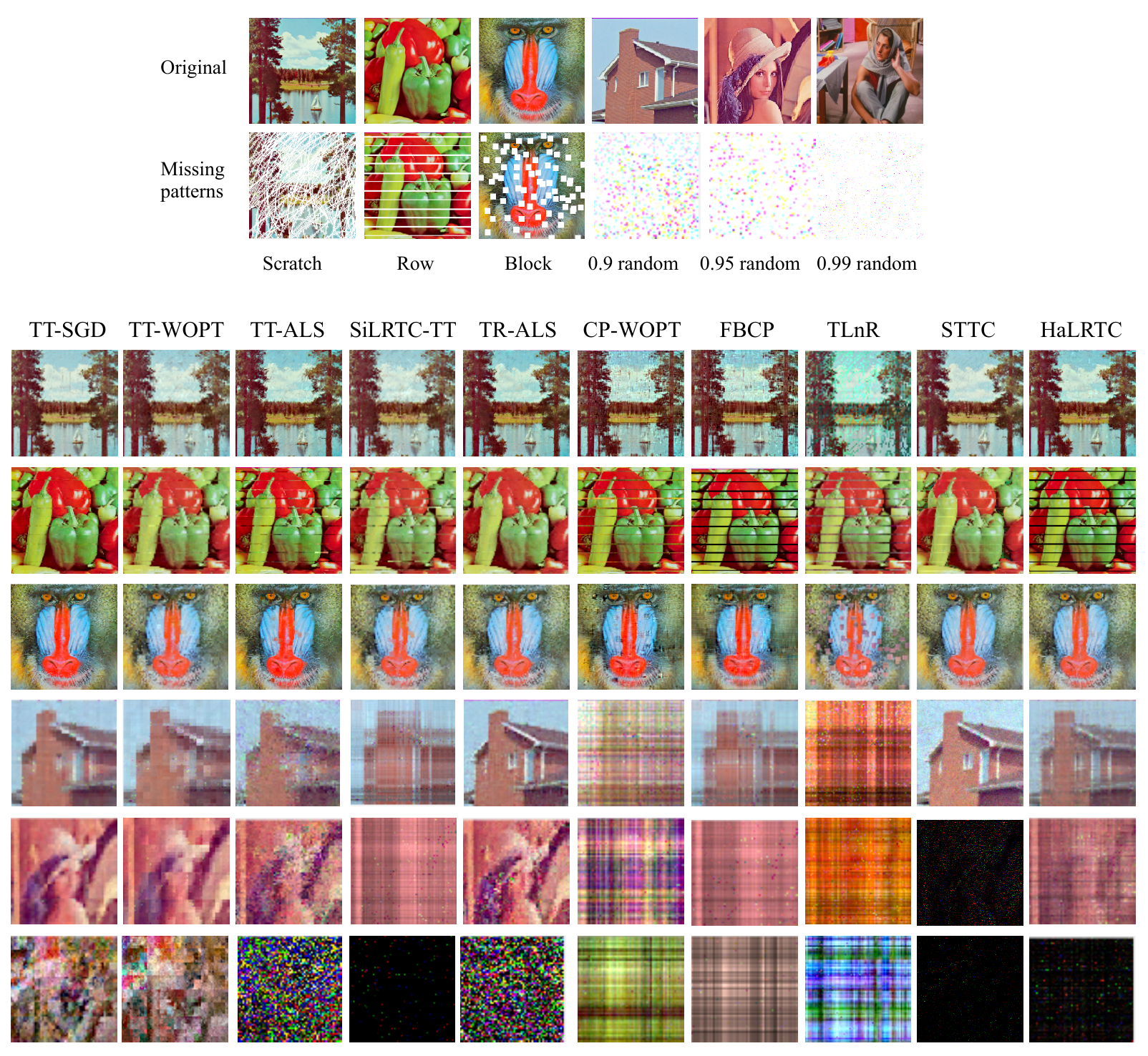}
\caption{The first and second row of the figure is the fully observed benchmark images and the corresponding missing patterns respectively, below which the visual completion results of the ten algorithms under the different missing patterns (i.e., scratch missing, row missing, block missing, 0.9 random missing, 0.95 random missing, and 0.99 random missing) are shown.}
\label{gt}
\end{center}
\end{figure}

In this experiment, we consider several irregular missing cases (the scratch missing, the whole row missing and the block missing) and some high-random-missing cases on benchmark RGB images. The parameter settings for each compared algorithms are tuned to get the best performance. The completion results from Figure \ref{gt} and Table \ref{gtt} we can see, our algorithms show high completion performance in all the missing cases. Moreover, for irregular missing cases and 0.8 random missing cases, STTC and HaLRTC performs well and achieve low RSE values. However, the two algorithms fail to solve the completion task when the random missing rate is 0.9 and 0.99,  this is because the nuclear-norm-based and total-variations-based algorithms cannot explore low-rank and local information when only a very small amount of entries is obtained. It should be noted that the 0.99 random missing case is a challenging task among all the image completion algorithms. Our two proposed algorithms with VDT method can achieve high performance under this situation while the other algorithms fail.

\begin{table}[h]
\centering
\caption{Comparison of the inpainting performance (RSE and PSNR) of ten algorithms under six missing situations.}
\label{gtt}
\resizebox{\textwidth}{25mm}{
\begin{tabular}{c|c|c|c|c|c|c|c|c|c|c|c}
\hline
\hline
 missing patterns&indices& TT-SGD & TT-WOPT &TT-ALS&SiLRTC-TT&TR-ALS& CP-WOPT & FBCP & TLnR&STTC&HaLRTC\\
 \hline
Scratch&\makecell[cc]{RSE\\PSNR}   &     \makecell[cc]{\textbf{0.09946}\\\textbf{24.93}}    &     \makecell[cc]{ 0.1455\\ 21.62}  &\makecell[cc]{0.1319 \\22.48 }       & \makecell[cc]{0.1522 \\21.23 }    &\makecell[cc]{ 0.1173\\23.49 }  &\makecell[cc]{ 0.2160\\18.19 } &\makecell[cc]{0.1185 \\23.40 } &\makecell[cc]{ 0.2871\\ 15.72} &\makecell[cc]{ 0.1085\\ 24.17} &\makecell[cc]{ 0.1168\\23.53 }  \\
\hline
Row &\makecell[cc]{RSE\\PSNR}&\makecell[cc]{ \textbf{0.07319}\\ \textbf{27.91}} &\makecell[cc]{0.09629 \\ 25.54} &\makecell[cc]{0.1385 \\ 22.38} &\makecell[cc]{0.1379 \\ 22.41}  &  \makecell[cc]{ 0.07325\\ 27.91} &\makecell[cc]{ 0.1653\\20.84 } &\makecell[cc]{0.3605 \\14.07 } &\makecell[cc]{ 0.1797\\20.12 } &\makecell[cc]{0.1069 \\ 24.62} &\makecell[cc]{ 0.3605\\ 14.07} \\
\hline
Block  &\makecell[cc]{RSE\\PSNR}&\makecell[cc]{0.08084\\27.20} &\makecell[cc]{0.09196 \\ 26.09} &\makecell[cc]{ 0.09671\\25.65 } &\makecell[cc]{0.09511 \\29.86 } &\makecell[cc]{ 0.08517\\26.75 } &\makecell[cc]{0.1391 \\22.49 } &\makecell[cc]{0.1147 \\ 24.17} &\makecell[cc]{0.1579 \\21.39 } &\makecell[cc]{\textbf{0.07315} \\\textbf{28.07} }&\makecell[cc]{ 0.08167\\27.12 }  \\
\hline
%0.8 random &\makecell[cc]{RSE\\PSNR}&\makecell[cc]{\textbf{ 0.09324}\\\textbf{22.50} }  &\makecell[cc]{0.1005 \\ 21.85} &\makecell[cc]{0.1351 \\19.28 } &\makecell[cc]{0.1362 \\19.21 } &\makecell[cc]{0.1388 \\19.05 }&\makecell[cc]{0.2564 \\13.72 } &\makecell[cc]{ 0.1321\\ 19.48} &\makecell[cc]{ 0.8625\\3.17 } &\makecell[cc]{\textbf{0.09514} \\\textbf{22.33 }} &\makecell[cc]{0.1485 \\ 18.46} \\
%\hline
0.9 random &\makecell[cc]{RSE\\PSNR}&\makecell[cc]{0.1444\\21.11 } &\makecell[cc]{0.1635 \\20.03 } &\makecell[cc]{ 0.1891\\18.77 } &\makecell[cc]{ 0.1969\\18.41 }   & \makecell[cc]{ \textbf{0.1090}\\ \textbf{23.55}} &\makecell[cc]{0.3209 \\14.17 } &\makecell[cc]{0.1967 \\18.43 } &\makecell[cc]{ 0.5815\\ 9.01} &\makecell[cc]{ 0.1845\\18.98 } &\makecell[cc]{ 0.1621\\ 20.11} \\
\hline
0.95 random &\makecell[cc]{RSE\\PSNR}&\makecell[cc]{\textbf{0.1576} \\\textbf{21.17} } &\makecell[cc]{0.1797\\20.32} &\makecell[cc]{0.2547 \\17.00 } &\makecell[cc]{0.2865 \\15.98 }  & \makecell[cc]{ 0.2147\\18.49 } &\makecell[cc]{ 0.4045\\12.98 } &\makecell[cc]{ 0.2850\\ 16.02} &\makecell[cc]{ 0.5557\\10.23 } &\makecell[cc]{- \\ -} &\makecell[cc]{0.2820 \\16.11 } \\
\hline
0.99 random &\makecell[cc]{RSE\\PSNR}&\makecell[cc]{0.3318 \\15.81} &\makecell[cc]{ \textbf{0.2520}\\\textbf{15.30} } &\makecell[cc]{- \\- } &\makecell[cc]{0.4049 \\6.98}  &\makecell[cc]{- \\ -} &\makecell[cc]{ 0.4749\\12.70 } &\makecell[cc]{0.4074 \\ 14.03} &\makecell[cc]{0.8545 \\ 7.30} &\makecell[cc]{ -\\- } &\makecell[cc]{0.9129 \\ 7.03} \\
\hline\hline
\end{tabular}
}
\end{table}

\subsection{Video and Hyperspectral Image Completion}
For large-scale data completion task, we test a video and a hyperspectral image (HSI) in the following experiments. For our proposed algorithms, we only test TT-SGD because TT-SGD is better for large-scale data than TT-WOPT. In addition, when large-scale data is employed, many algorithms which work well on benchmark images will become inefficient or ineffective, so we compare TT-SGD to only several algorithms (TT-ALS, CP-WOPT, FBCP, and HaLRTC).

\begin{figure}[htp]
\begin{center}
\includegraphics[width=1\linewidth]{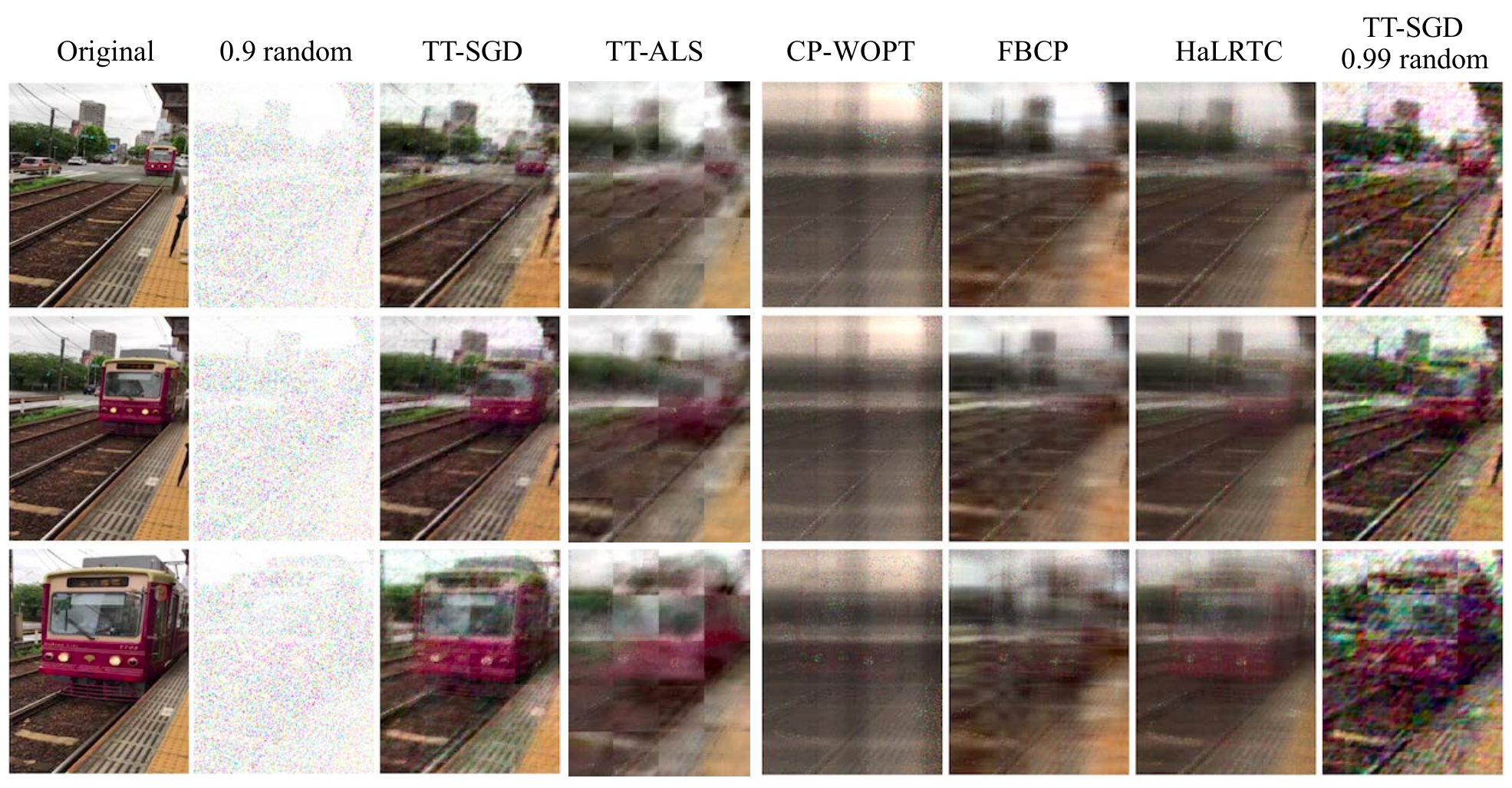}
\caption{Video completion results of TT-SGD, TT-ALS, CP-WOPT, FBCP, and HaLRTC under random missing cases. The first row to the last row show the completion results of the 1st frame, the 75th frame and the 100th frame of the video respectively.}
\label{video1}
\end{center}
\end{figure}

First, we test a video which records a moving train. The size of the data is $320\times 256 \times 3 \times 100$ and the background of the video changes by frames. By the VDT method, we first reshape the data to size $2 \times2\times 2\times 2\times 2\times 2\times 5\times 2\times 2\times 2\times 2\times 2\times 2\times 4\times 3\times 100$, then permute it by index $\{1\ 8\ 2\ 9\ 3\ 10\ 4\ 11\ 5\ 12\ 6\ 13\ 7\ 14\ 15\ 16\}$, and finally we reshape it to size $4\times4\times4\times4\times4\times4\times20\times3\times100$ as the input tensor. We compare three random missing cases ($m_r=0.7$, $m_r=0.9$ and $m_r=0.99$) in this experiment. Part of the visual results are shown in Figure \ref{video1}, and the numerical results are shown in Table \ref{vm}. The performance of TT-SGD outperforms other compared algorithms. More specifically, it can recover the video well even there is only 1\% sampled entries while other compared algorithms fail in this high missing rate case. It should also be noted that the time cost of TT-SGD is lower than the other compared algorithms, which shows high efficiency of TT-SGD.

\begin{table}[htp]\footnotesize
\centering
\caption{Numerical results (RSE and PSNR) on video completion experiments of five algorithms under three random missing cases. }
\label{vm}
\begin{tabular}{c|c|c|c}
\hline \hline
& $m_r=0.7$ &$m_r=0.9$ &$m_r=0.99$  \\
\hline
Algorithm&RSE  PSNR Time &RSE  PSNR Time &RSE  PSNR Time \\
\hline
TT-SGD &\textbf{0.1459 22.67 680.94}&\textbf{0.2045 19.87 674.17}&\textbf{0.2185 19.24 698.11}   \\
\hline
TT-ALS & 0.2116 19.48 7100.43&0.2400 18.39  1622.42&0.2557 17.8466  793.76\\
\hline
CP-WOPT &0.2673 17.41 825.06& 0.3264 15.67 790.60&0.3610 14.80 814.44\\
\hline
FBCP & 0.2204 19.11 870.89&0.2547 17.86 920.78&0.3258 15.72 720.01\\
\hline
HaLRTC& 0.1758 21.16 1132.05&  0.2562 17.78 1044.88&0.8844 7.016 1121.37\\
\hline \hline
\end{tabular}
\end{table}

\begin{figure}[H]
\begin{center}
\includegraphics[width=0.8\linewidth]{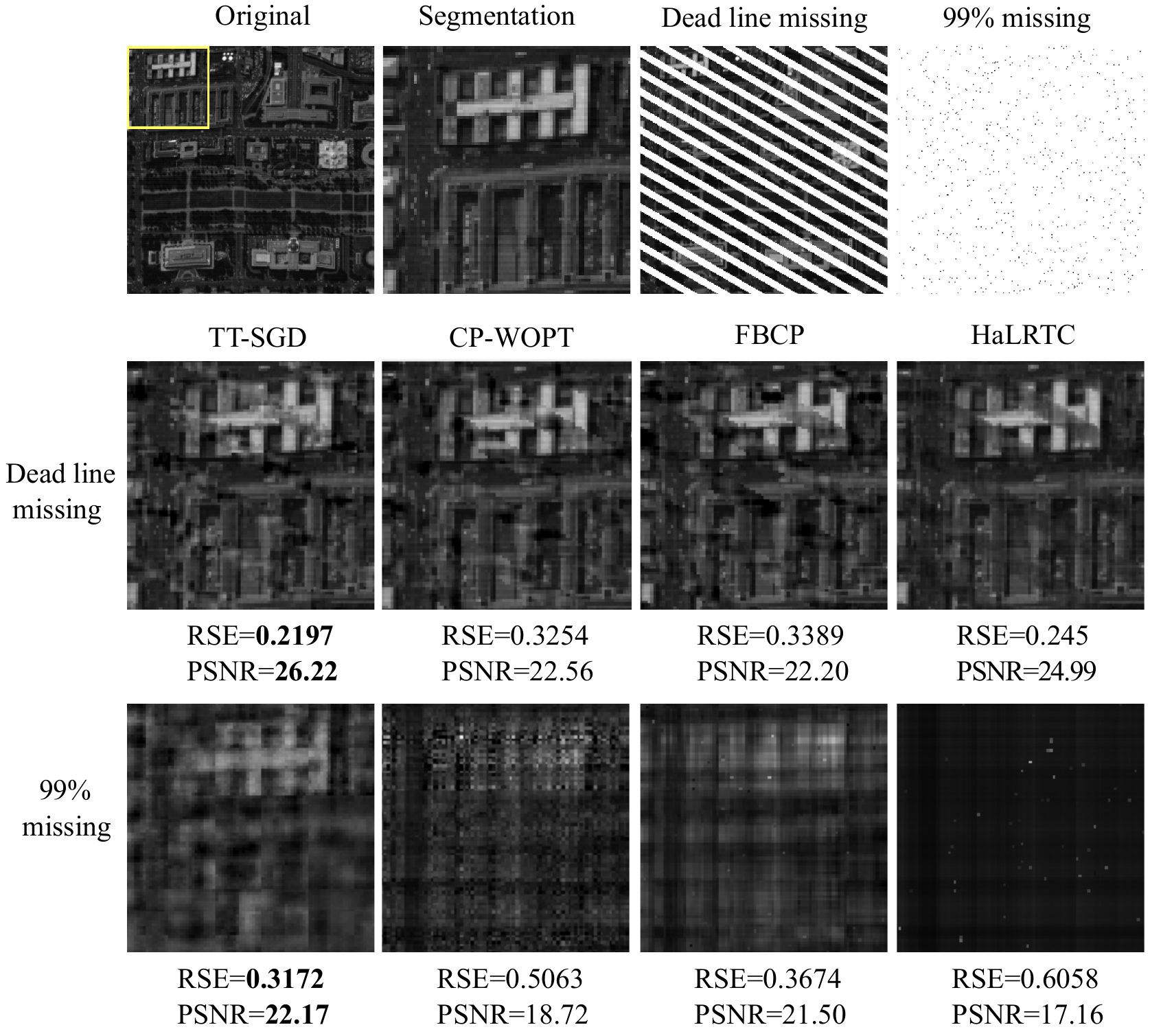}
\caption{HSI completion results of the four algorithms. We show the image of the first channel of the HSI. The first row is the original image, the segmentation to show the completion performance, the dead line missing pattern, and the 0.99 random missing pattern. The second row and the third row show the completion results.}
\label{hsi}
\end{center}
\end{figure}

Then we test TT-SGD, CP-WOPT, FBCP and HaLRTC on a hyperspectral image (HSI) of size $256\times 256 \times191$ recorded by a satellite. Due to the inferior working condition of satellite sensors, the collected data often has Gaussian noise, impulse noise, dead lines, and stripes \cite{zhang2014hyperspectral}. In this experiment, we first consider the situation when the HSI has `dead lines', which is a common missing case in HSI record. Then we consider the case when only 1\% of the data is obtained, which is meaningful in data compression and transformation. We transform the HSI data to $16\times16 \times16 \times16 \times191$ by VDT method as the input for TT-SGD and apply original three-order tensor as the input for the other compared algorithms. We set TT-ranks as $48$ and $24$ for dead line missing case and 99\% missing case respectively. The visual completion results in Figure \ref{hsi} shows the image of the first channel of the HSI and the numerical results are the evaluation of the overall completion performance.

TT-SGD performs best among the algorithms at both dead line missing case and 99\% random missing case. In 99\% random missing case, HaLRTC fails the completion task,  while CP-WOPT and FBCP obtain lower performance than TT-SGD. In addition, it should be noted that the volume of data is about $1.25\times 10^7$, and when the iteration reaches $1\times10^6$ (16\% of the total data), the optimization of TT-SGD is converged. This indicates that TT-SGD has fast and efficient computation.

\section{Conclusion}
In this paper, in order to solve the tensor completion problem, based on tensor train decomposition and gradient descent method, we propose two tensor completion algorithms named TT-WOPT and TT-SGD. We first cast the completion problem into solving the optimization models, then we use gradient descent methods to find the optimal core tensors of TT decomposition. Finally, the TT core tensors are applied to approximate the missing entries of the incomplete tensor. Furthermore, to improve the performance of the proposed algorithms, we propose the VDT method to tensorize visual data to higher-order. We conduct simulation experiments and visual data experiments to compare our algorithms to the state-of-the-art algorithms. From the simulation experiments we can see, the performance of our algorithms stays stable when the tensor order increases. Moreover, the visual data experiments show that after higher-order tensorization by VDT, the performance of our two algorithms can be improved. Our algorithms outperform the compared state-of-the-art algorithms in various missing situations, particularly when the tensor order is high and the missing rate is high. More specially, our algorithms with VDT method can process extreme high random missing situation (i.e., 99\% random missing) well while other algorithms fail. Besides, our proposed TT-SGD achieves low computational complexity and high efficiency in processing large-scale data.

The high performance of the proposed algorithms shows that TT-based tensor completion is a promising aspect. It should be noted that TT-rank setting is essential to obtain better experiment results and it is selected manually in common. We will extend our algorithms by choosing TT-rank automatically in our future work.
 
\section*{Acknowledgement}
This work was supported by JSPS KAKENHI (Grant No. 17K00326, 15H04002, 18K04178), JST CREST (Grant No. JPMJCR1784) and the National Natural Science Foundation of China (Grant No. 61773129).
\section*{References}

\bibliography{SPIC}

\end{document}